\documentclass[12pt]{amsart}
\usepackage{amssymb,amsmath,amsfonts,latexsym,setspace}
\usepackage{bm}
\newcommand{\newsection}[1]{\setcounter{equation}{0} \section{#1}}
\newcommand{\bea}{\begin{eqnarray}}
\newcommand{\eea}{\end{eqnarray}}

\newcommand{\cle}{\mathcal{E}}

\newcommand{\clh}{\mathcal{H}}

\newcommand{\cll}{\mathcal{L}}
\newcommand{\clm}{\mathcal{M}}

\newcommand{\clq}{\mathcal{Q}}
\newcommand{\clr}{\mathcal{R}}
\newcommand{\cls}{\mathcal{S}}

\newcommand{\raro}{\rightarrow}

\def \qed {\hfill \vrule height6pt width 6pt depth 0pt}
\def\textmatrix#1&#2\\#3&#4\\{\bigl({#1 \atop #3}\ {#2 \atop #4}\bigr)}
\def\dispmatrix#1&#2\\#3&#4\\{\left({#1 \atop #3}\ {#2 \atop #4}\right)}
\newcommand{\be}{\begin{equation}}
\newcommand{\ee}{\end{equation}}
\newcommand{\ben}{\begin{eqnarray*}}
\newcommand{\een}{\end{eqnarray*}}

\newcommand{\NI}{\noindent}

\newcommand{\bi}{\begin{itemize}}
\newcommand{\ei}{\end{itemize}}

\newtheorem{Theorem}{\sc Theorem}[section]
\newtheorem{Lemma}[Theorem]{\sc Lemma}
\newtheorem{Proposition}[Theorem]{\sc Proposition}
\newtheorem{Corollary}[Theorem]{\sc Corollary}
\newtheorem{Definition}[Theorem]{\sc Definition}
\newtheorem{Example}[Theorem]{\sc Example}
\newtheorem{Remark}[Theorem]{\sc Remark}
\newtheorem{Note}[Theorem]{\sc Note}
\newtheorem{Question}{\sc Question}
\newtheorem{ass}[Theorem]{\sc Assumption}
\newcommand{\bt}{\begin{Theorem}}
\def\beginlem{\begin{Lemma}}
\def\beginprop{\begin{Proposition}}
\def\begincor{\begin{Corollary}}
\def\begindef{\begin{Definition}}
\def\beginexamp{\begin{Example}}
\def\beginrem{\begin{Remark}}
\def\beginq{\begin{Question}}
\def\beginass{\begin{ass}}
\def\beginnote{\begin{Note}}
\newcommand{\et}{\end{Theorem}}
\def\endlem{\end{Lemma}}
\def\endprop{\end{Proposition}}
\def\endcor{\end{Corollary}}
\def\enddef{\end{Definition}}
\def\endexamp{\end{Example}}
\def\endrem{\end{Remark}}
\def\endq{\end{Question}}
\def\endass{\end{ass}}
\def\endnote{\end{Note}}

\begin{document}

\title[Jordan Blocks of $H^2(\mathbb{D}^n)$]{Jordan Blocks of $H^2(\mathbb{D}^n)$}

\author[Jaydeb Sarkar]{Jaydeb Sarkar}
\address{Indian Statistical Institute, Statistics and Mathematics Unit, 8th Mile, Mysore Road, Bangalore, 560059, India}
\email{jay@isibang.ac.in, jaydeb@gmail.com}




\begin{abstract}
We develop a several variables analog of the Jordan blocks of the
Hardy space $H^2(\mathbb{D})$. In this consideration, we obtain a
complete characterization of the doubly commuting quotient modules
of the Hardy module $H^2(\mathbb{D}^n)$. We prove that a quotient
module $\clq$ of $H^2(\mathbb{D}^n)$ ($n \geq 2$) is doubly
commuting if and only if \[\clq = \clq_{\Theta_1} \otimes \cdots
\otimes \clq_{\Theta_n},\]where each $\clq_{\Theta_i}$ is either a
one variable Jordan block $H^2(\mathbb{D})/\Theta_i H^2(\mathbb{D})$
for some inner function $\Theta_i$ or the Hardy module
$H^2(\mathbb{D})$ on the unit disk for all $i = 1, \ldots, n$. We
say that a submodule $\cls$ of $H^2(\mathbb{D}^n)$ is co-doubly
commuting if the quotient module $H^2(\mathbb{D}^n)/\cls$ is doubly
commuting. We obtain a Beurling like theorem for the class of
co-doubly commuting submodules of $H^2(\mathbb{D}^n)$. We prove that
a submodule $\cls$ of $H^2(\mathbb{D}^n)$ is co-doubly commuting if
and only if
\[\cls = \mathop{\sum}_{i=1}^m \Theta_i
H^2(\mathbb{D}^n),\]for some integer $m \leq n$ and one variable
inner functions $\{\Theta_i\}_{i=1}^m$.
\end{abstract}

\maketitle

\newsection{Introduction}

Let $\mathbb{C}^n = \mathbb{C} \times \cdots \times \mathbb{C}$ be
the $n$-dimensional complex Euclidean space with $n \geq 1$ and
$\mathbb{D}^n = \{(z_1, \ldots, z_n) : |z_i| < 1, i = 1, \ldots,
n\}$ be the open unit polydisc. We denote the elements of
$\mathbb{C}^n$ by $\bm{z} = (z_1, \ldots, z_n)$ where $z_i \in
\mathbb{C}$ for all $i= 1, \ldots, n$. The \textit{Hardy space}
$H^2(\mathbb{D}^n)$ on the polydisc is the Hilbert space of all
holomorphic functions $f$ on $\mathbb{D}^n$ such that
\[\|f\|_{H^2(\mathbb{D}^n)} := \bigg(\mathop{\mbox{sup}}_{0 \leq r < 1}
\mathop{\int}_{\mathbb{T}^n} |f(r \bm{z})|^2
d\bm{\theta}\bigg)^{\frac{1}{2}} < \infty,\]where $d\bm{\theta}$ is
the normalized Lebesgue measure on the torus $\mathbb{T}^n$, the
distinguished boundary of $\mathbb{D}^n$ and $r\bm{z} : = (rz_1,
\ldots, r z_n)$ (cf. \cite{R}, \cite{FS}).

The multiplication operators by the coordinate functions turns
$H^2(\mathbb{D}^n)$ into a \textit{Hilbert module} over
$\mathbb{C}[\bm{z}] = \mathbb{C}[z_1, \ldots, z_n]$, the ring of
polynomials in $n$ variables with complex coefficients in the
following sense:
\[\mathbb{C}[\bm{z}] \times H^2(\mathbb{D}^n) \raro
H^2(\mathbb{D}^n),\quad (p, f) \mapsto p(M_{z_1}, \ldots,
M_{z_n})f,\]for all $p \in \mathbb{C}[\bm{z}]$ and $f \in
H^2(\mathbb{D}^n)$ (cf. \cite{DP}). We also call the Hilbert module
$H^2(\mathbb{D}^n)$ over $\mathbb{C}[\bm{z}]$ as the \textit{Hardy
module}. A closed subspace $\cls \subseteq H^2(\mathbb{D}^n)$ is
said to be a \textit{submodule} of $H^2(\mathbb{D}^n)$ if $M_{z_i}
\cls \subseteq \cls$ for all $i = 1, \ldots, n$. A closed subspace
$\clq \subseteq H^2(\mathbb{D}^n)$ is said to be a \textit{quotient
module} of $H^2(\mathbb{D}^n)$ if $\clq^{\perp} ( \cong
H^2(\mathbb{D}^n)/\clq )$ is a submodule of $H^2(\mathbb{D}^n)$.

Let $\cls$ be a submodule and $\clq$ be a quotient module of
$H^2(\mathbb{D}^n)$. Then the module multiplication operators on
$\cls$ and $\clq$ are given by the restrictions $(R_{z_1}, \ldots,
R_{z_n})$ and the compressions $(C_{z_1}, \ldots, C_{z_n})$ of the
module multiplications of $H^2(\mathbb{D}^n)$, respectively. That
is,
\[R_{z_i} = M_{z_i}|_{\cls} \quad \mbox{and} \quad C_{z_i} = P_{\clq}
M_{z_i}|_{\clq},\]for all $i = 1, \ldots, n$. Here, for a given
closed subspace $\clm$ of a Hilbert space $\clh$, we denote the
orthogonal projection of $\clh$ onto $\clm$ by $P_{\clm}$. Note that
\[R_{z_i}^* = P_{\cls} M_{z_i}^*|_{\cls} \quad \quad \mbox{and}
\quad \quad C_{z_i}^* = M_{z_i}^*|_{\clq},\]for all $i= 1, \ldots,
n$.

\vspace{0.2in}

\NI\textsf{Jordan blocks of $H^2(\mathbb{D})$:} A closed subspace
$\clq \subseteq H^2(\mathbb{D})$ is said to be a \textit{Jordan
block} of $H^2(\mathbb{D})$ if $\clq$ is a quotient module and $\clq
\neq H^2(\mathbb{D})$ (see \cite{NF1}, \cite{NF}). By Beurling's
theorem \cite{B}, a closed subspace $\clq (\neq H^2(\mathbb{D}))$ is
a quotient module of $H^2(\mathbb{D})$ if and only if the submodule
$\clq^{\perp}$ is given by $\clq^{\perp} = \Theta H^2(\mathbb{D})$
for some inner function $\Theta \in H^{\infty}(\mathbb{D})$. In
other words, the quotient modules and hence the Jordan blocks of
$H^2(\mathbb{D})$ are precisely given by
\[\clq_\Theta : = H^2(\mathbb{D})/\Theta H^2(\mathbb{D}),\]for inner
functions $\Theta \in H^{\infty}(\mathbb{D})$.

\vspace{0.2in}

\NI\textsf{Jordan blocks of $H^2(\mathbb{D}^n)$ ($n >1)$:} First we
note that the Hardy module $H^2(\mathbb{D}^n)$ (with $n
>1$) can be identified with the $n$-fold Hilbert space tensor product of
the Hardy space $H^2(\mathbb{D})$ on the disc
\[ \underbrace{H^2(\mathbb{D}) \otimes \cdots \otimes
H^2(\mathbb{D})}_{n \rm \; times},\]via the unitary map $U :
H^2(\mathbb{D}^n) \raro H^2(\mathbb{D}) \otimes \cdots \otimes
H^2(\mathbb{D})$, where $U(z_1^{l_1} \cdots z_n^{l_n}) := z^{l_1}
\otimes \cdots \otimes z^{l_n}$ for all $l_1, \ldots, l_n \in
\mathbb{N}$. Moreover, $H^2(\mathbb{D}) \otimes \cdots \otimes
H^2(\mathbb{D})$ is a Hilbert module over $\mathbb{C}[\bm{z}]$ with
the module multiplication operators \[\{I_{{H^2(\mathbb{D})}}
\otimes \cdots \otimes \underbrace{M_z}_{i^{\rm th}}  \otimes \cdots
\otimes I_{{H^2(\mathbb{D})}}\}_{i=1}^n.\] Therefore, that $U$ is a
module map \[U M_{z_i} = (I_{{H^2(\mathbb{D})}} \otimes \cdots
\otimes \underbrace{M_z}_{i^{\rm th}} \otimes \cdots \otimes
I_{{H^2(\mathbb{D})}}) U,\] for all $1 \leq i \leq n$. It is easy to
see that
\[M_{z_i} M_{z_j}^* = M_{z_j}^* M_{z_i},\]for all $i \neq j$.

\NI The above fact is one of the motivations to introduce the
following notion and the title of the paper.

\begin{Definition} Let $\clq$ be a quotient module of $H^2(\mathbb{D}^n)$
and $n >1$. Then $\clq$ is said to be a Jordan block of
$H^2(\mathbb{D}^n)$ if $\clq$ is doubly commuting, that is, $C_{z_i}
C_{z_j}^* = C_{z_j}^* C_{z_i}$, for all $1 \leq i < j \leq n$ and
$\clq \neq H^2(\mathbb{D}^n)$. Also a closed subspace $\cls$ of
$H^2(\mathbb{D}^n)$ is said to be co-doubly commuting submodule of
$H^2(\mathbb{D}^n)$ if $H^2(\mathbb{D}^n)/\cls$ is a doubly
commuting quotient module.
\end{Definition}

\NI However, in most of the following we will simply regard a Jordan
block of $H^2(\mathbb{D}^n)$ as a doubly commuting quotient module
of $H^2(\mathbb{D}^n)$.

The study of the doubly commuting quotient modules of the Hardy
module $H^2(\mathbb{D}^2)$ was initiated by Douglas and Yang in
\cite{DY1} and \cite{DY} (also see \cite{BCL}). Later in \cite{INS}
Izuchi, Nakazi and Seto obtained a classification of the doubly
commuting quotient modules of the Hardy module $H^2(\mathbb{D}^2)$
(see Theorems 2.1 and 3.1 in \cite{INS}).

In this paper we completely classify the doubly commuting quotient
modules of $H^2(\mathbb{D}^n)$ for any $n \geq 2$. In this
consideration, we provide a more refined analysis compared to
\cite{INS}. More specifically, our method is based on the Hilbert
tensor product structure of the Hardy module $H^2(\mathbb{D}^n)$
which also yield new proofs of earlier results by Izuchi, Nakazi and
Seto \cite{INS} concerning the base case $n = 2$.

A key example of doubly commuting quotient modules  over
$\mathbb{C}[\bm{z}]$ is the Hilbert tensor product of $n$ quotient
modules of the Hardy module $H^2(\mathbb{D})$. That is, if we
consider $n$ quotient modules $\{\clq_i\}_{i=1}^n$ of the Hardy
module $H^2(\mathbb{D})$ then
\[\clq = \clq_1 \otimes \cdots \otimes \clq_n,\] is a doubly
commuting quotient module of $H^2(\mathbb{D}^n)$ with the module
multiplication operators as
\[\{I_{\clq_1} \otimes \cdots \otimes
\underbrace{P_{\clq_i} M_{z}|_{\clq_i}}_{i^{\rm\, th} \rm\,place}
\otimes \cdots \otimes I_{\clq_n}\}_{i=1}^n.\] We prove that a
doubly commuting quotient module of $H^2(\mathbb{D}^n)$ can be also
represented by the Hilbert space tensor product of quotient modules
of $H^2(\mathbb{D})$ in the above form. This result is then used to
prove a Beurling type theorem for the co-doubly commuting
submodules.

We now summarize the contents of this paper. In Section 2, we give
relevant background for the main results of this paper. In Section
3, we prove that a quotient module $\clq$ of $H^2(\mathbb{D}^n)$ is
doubly commuting if and only if $\clq$ is the $n$ times Hilbert
tensor product of quotient modules of the Hardy module
$H^2(\mathbb{D})$. In Section 4, we characterize the class of
co-doubly commuting submodules of $H^2(\mathbb{D}^n)$.

\newsection{Preparatory results }

In this section, we gather together some concepts and results
concerning various aspects of the Hardy modules that are used
frequently in the rest of this paper. Some of the results of the
present section are of independent interest.

We first recall that the module multiplication of the Hardy module
$H^2(\mathbb{D})$ satisfies the following relation \[M_z M_z^* =
I_{H^2(\mathbb{D})} - P_{\mathbb{C}},\]where $P_{\mathbb{C}}$
denotes the orthogonal projection of $H^2(\mathbb{D})$ onto the
space of constant functions. Moreover, if $\clq_{\Theta} =
H^2(\mathbb{D})/\Theta H^2(\mathbb{D})$ is a Jordan block for some
inner function $\Theta \in H^\infty(\mathbb{D})$, then we have
\[P_{\clq_{\Theta}} = I_{H^2(\mathbb{D})} - M_{\Theta} M^*_{\Theta}
\quad \; \mbox{and} \; \quad P_{\Theta H^2(\mathbb{D})} = M_{\Theta}
M_{\Theta}^*.\]We also have \[I_{\clq} - C_z C_z^* =
P_{\clq}(I_{H^2(\mathbb{D})} - M_z M_z^*)|_{\clq} = P_{\clq}
P_{\mathbb{C}}|_{\clq},\]where $\clq$ is a quotient module of
$H^2(\mathbb{D})$.

The following lemma is well known.

\begin{Lemma}\label{J-Q}
Let $\clq_{\Theta}$ be a Jordan block of $H^2(\mathbb{D})$ for some
inner function $\Theta \in H^{\infty}(\mathbb{D})$. Then
\[P_{\clq_{\Theta}} 1 = 1 - \overline{\Theta(0)} \Theta,\]and
\[(P_{\clq_{\Theta}} P_{\mathbb{C}} P_{\clq_{\Theta}}) 1 = (1 - |\Theta(0)|^2) (1 - \overline{\Theta(0)}
\Theta).\]
\end{Lemma}

\NI\textsf{Proof.} By virtue of $M_{\Theta}^* 1 =
\overline{\Theta(0)}$ we have
\[P_{\clq_{\Theta}} 1 = ( I_{H^2(\mathbb{D})} - M_{\Theta} M_{\Theta}^*) 1 = 1 - M_{\Theta} (M_{\Theta}^* 1) = 1 - \overline{\Theta(0)}
\Theta.\]For the second equality, we compute \[(P_{\clq_{\Theta}}
P_{\mathbb{C}} P_{\clq_{\Theta}}) 1 = (P_{\clq_{\Theta}}
P_{\mathbb{C}}) (1 - \overline{\Theta(0)} \Theta) =
P_{\clq_{\Theta}} (1 - |\Theta(0)|^2) = (1 - |\Theta(0)|^2) (1 -
\overline{\Theta(0)} \Theta).\] This completes the proof. \qed

This lemma has the following immediate corollary.

\begin{Corollary}\label{Q_1}
Let $\clq$ be a quotient module of $H^2(\mathbb{D})$. Then
\[P_{\clq} 1 \in \mbox{ran}
(P_{\clq} P_{\mathbb{C}} P_{\clq}).\]
\end{Corollary}
\NI \textsf{Proof.} If $\clq = H^2(\mathbb{D})$ then the result
follows trivially. If $\clq \neq H^2(\mathbb{D})$ then $\clq$ is a
Jordan block and hence the conclusion follows from Lemma \ref{J-Q}.
\qed

The following lemma is a variation on the theme of the isometric
dilation theory of contractions.

\begin{Lemma}\label{Q-L}
Let $\clq$ be a quotient module of $H^2(\mathbb{D})$ and $\cll =
\mbox{ran} (I_{\clq} - C_z C^*_z) = \mbox{ran} (P_{\clq}
P_{\mathbb{C}} P_{\clq})$. Then
\[\clq = \mathop{\vee}_{l=0}^{\infty} P_{\clq} M_z^l \cll.\]
\end{Lemma}
\NI\textsf{Proof.} The result is trivial if $\clq = \{0\}$. Let
$\clq \neq \{0\}$. Notice that
\[\mathop{\vee}_{l=0}^{\infty} P_{\clq} M_z^l \cll \subseteq \clq.\]
Let now $f \in \clq$ be such that $f \perp \vee_{l=0}^{\infty}
P_{\clq} M_z^l \cll$. Then for all $l \geq 0$ we have that $f \perp
P_{\clq} M_z^l P_{\clq} P_{\mathbb{C}} \clq$, or equivalently,
$P_{\mathbb{C}} M_z^{*l} f \in {\clq}^{\perp}$. Since
${\clq}^{\perp}$ is a proper submodule of $H^2(\mathbb{D})$, it
follows that
\[P_{\mathbb{C}} M_z^{*l} f = 0,\]for all $l \geq 0$. Consequently,
\[f = 0.\]This concludes the proof. \qed

In the following, we employ the standard multi-index notation that
$\mathbb{N}^n = \{(k_1, \ldots, k_n) : k_i \in \mathbb{N}, i = 1,
\ldots, n\}$ and for any $\bm{k} = (k_1, \ldots, k_n) \in
\mathbb{N}^n$ we denote $\bm{z}^{\bm{k}} = z_1^{k_1} \cdots
z_n^{k_n}$ and $M_{\bm{z}}^{\bm{k}} = M_{z_1}^{k_1} \cdots
M_{z_n}^{k_n}$.

We now present a characterization of $M_{z_1}$-reducing subspace of
$H^2(\mathbb{D}^n)$. However the technique used here seems to be
well known in the study of the reducing subspaces.

\begin{Proposition}\label{reducing}
Let $n >1$ and $\cls$ be a closed subspace of $H^2(\mathbb{D}^n)$.
Then $\cls$ is a $(M_{z_2}, \ldots, M_n)$-reducing subspace of
$H^2(\mathbb{D}^n)$ if and only if
\[\cls = \cls_1 \otimes \underbrace{H^2(\mathbb{D}) \otimes \cdots \otimes H^2(\mathbb{D})}_{(n-1) \rm \,times},\]for some closed subspace $\cls_1$ of
$H^2(\mathbb{D})$.
\end{Proposition}

\NI\textsf{Proof.} Let $\cls$ be a $(M_{z_2}, \ldots,
M_{z_{n}})$-reducing closed subspace of $H^2(\mathbb{D}^n)$, that
is, for all $2 \leq i \leq n$ we have
\[M_{z_i} P_{\cls} = P_{\cls} M_{z_i}.\] Following Agler's
hereditary functional calculus (cf. \cite{At})
\[\begin{split}(\mathop{\Pi}_{i=2}^{n} (1 - z_i
\bar{w}_i))(\bm{M_z}, \bm{M_z}) & = \sum_{\substack{
0 \leq i_1 < \ldots < i_l \leq n\\
i_1, i_2 \neq 1}} (-1)^l (z_{i_1} \cdots z_{i_l} \bar{w}_{i_1}
\cdots \bar{w}_{i_l}) (\bm{M_z}, \bm{M_z})\\ & = \sum_{\substack{
0 \leq i_1 < \ldots < i_l \leq n\\
i_1, i_2 \neq 1}} (-1)^l M_{z_{i_1}} \cdots M_{z_{i_l}}
M^*_{{z}_{i_1}} \cdots M^*_{{z}_{i_l}}\\ & = P_{H^2(\mathbb{D})}
\otimes P_{\mathbb{C}} \otimes \cdots \otimes P_{\mathbb{C}},
\end{split}\]where $\bm{M_z} = (M_{z_1}, \ldots, M_{z_n})$. Consequently, \[(P_{H^2(\mathbb{D})} \otimes P_{\mathbb{C}} \otimes
\cdots \otimes P_{\mathbb{C}}) P_{\cls} = P_{\cls}
(P_{H^2(\mathbb{D})} \otimes P_{\mathbb{C}} \otimes \cdots \otimes
P_{\mathbb{C}}),\]which yields that $P_{\cls} (P_{H^2(\mathbb{D})}
\otimes P_{\mathbb{C}} \otimes \cdots \otimes P_{\mathbb{C}})$ is an
orthogonal projection and \[P_{\cls} (P_{H^2(\mathbb{D})} \otimes
P_{\mathbb{C}} \otimes \cdots \otimes P_{\mathbb{C}}) =
(P_{H^2(\mathbb{D})} \otimes P_{\mathbb{C}} \otimes \cdots \otimes
P_{\mathbb{C}}) P_{\cls} = P_{\tilde{\cls_1}},\]where
$\tilde{\cls_1} := (H^2(\mathbb{D}) \otimes \mathbb{C} \otimes
\cdots \otimes \mathbb{C}) \cap \cls$. Let \[\tilde{\cls}_1 = \cls_1
\otimes \mathbb{C} \otimes \cdots \otimes \mathbb{C},\] for some
closed subspace $\cls_1$ of $H^2(\mathbb{D})$. We claim that
\[\cls = \overline{\mbox{span}} \{M_{z_2}^{l_2} \cdots M_{z_{n}}^{l_{n}}
\tilde{\cls}_1 : l_2, \ldots, l_{n} \in \mathbb{N}\} = \cls_1
\otimes H^2(\mathbb{D}^{n-1}).\] Now for any
\[f = \sum_{\bm{k} \in \mathbb{N}^n} a_{\bm{k}}
\bm{z}^{\bm{k}} \in \cls,\]we have \[f = P_{\cls} f = P_{\cls}
(\sum_{\bm{k} \in \mathbb{N}^n} M_{\bm{z}}^{\bm{k}} a_{\bm{k}}) =
\sum_{\bm{k} \in \mathbb{N}^n} M_{z_2}^{k_2}\cdots M_{z_{n}}^{k_{n}}
P_{\cls} (a_{\bm{k}} z_1^{k_1}),\]where $a_{\bm{k}} \in \mathbb{C}$
for all $\bm{k} \in \mathbb{N}^n$. But $P_{\cls} (a_{\bm{k}}
z_1^{k_1}) = P_{\cls} (P_{H^2(\mathbb{D})} \otimes P_{\mathbb{C}}
\otimes \cdots \otimes P_{\mathbb{C}}) (a_{\bm{k}} z_1^{k_1}) \in
\tilde{\cls}_1$ and hence $f \in \cls_1 \otimes
H^2(\mathbb{D}^{n-1})$. That is, $\cls \subseteq \cls_1 \otimes
H^2(\mathbb{D}^{n-1})$. On the other hand, since $\tilde{\cls}_1
\subseteq \cls$ and that $\cls$ is a $(M_{z_2}, \ldots,
M_{z_{n}})$-reducing subspace, we see that $\cls = \cls_1 \otimes
H^2(\mathbb{D}^{n-1})$.

\NI The converse part is immediate. This concludes the proof of the
proposition. \qed

The following result will be used in the final section.

\begin{Lemma}\label{P-F} Let $\{P_i\}_{i=1}^n$ be a collection of commuting orthogonal
projections on a Hilbert space $\clh$. Then \[\cll :=
\mathop{\sum}_{i=1}^n \mbox{ran} P_i,\] is closed and the orthogonal
projection of $\clh$ onto $\cll$ is given by
\[\begin{split}P_{\cll} & = P_1 (I - P_2) \cdots (I - P_n) +
P_2 (I - P_3) \cdots (I - P_n) + \cdots + P_{n-1} (I - P_n) + P_n\\
& = P_n (I - P_{n-1}) \cdots (I - P_1) + P_{n-1} (I - P_{n-2})
\cdots (I - P_1) + \cdots + P_2 (I - P_1) +
P_1.\end{split}\]Moreover,
\[P_{\cll}
= I - \mathop{\prod}_{i=1}^n (I - P_i).\]
\end{Lemma}
\NI\textsf{Proof.} We let \[O_i = P_i (I - P_{i+1}) \cdots (I -
P_{n-1}) (I - P_n),\] so that \[O_i = \mathop{\Pi}_{j= i+1}^n (I -
P_j) - \mathop{\Pi}_{j= i}^n (I - P_j),\]for all $i = 1, \ldots,
n-1$ and $O_n = P_n$. By the assumptions, $\{O_i\}_{i=1}^n$ is a
family of orthogonal projections with orthogonal ranges. We claim
that
\[\cll = \mbox{ran} O_1 \oplus \cdots \oplus
\mbox{ran} O_n.\]From the definition we see that $\cll \supseteq
\mbox{ran} O_1 \oplus \cdots \oplus \mbox{ran} O_n$. To prove the
reverse inclusion, first we observe that
\[\mathop{\sum}_{i=1}^n O_i = I - \mathop{\Pi}_{i=1}^n (I - P_i).\]
Now let $f = f_1 + \cdots + f_n \in \cll$ where $f_i \in \mbox{ran}
P_i$ for all $i=1,\ldots, n$. Then
\[\begin{split} (\mathop{\sum}_{i=1}^n O_i)f & = ( I - \mathop{\Pi}_{i=1}^n
(I - P_i)) f =  f - \mathop{\Pi}_{i=1}^n (I - P_i)f \\ & = f -
\mathop{\sum}_{j=1}^n \mathop{\Pi}_{i=1}^n (I - P_i)f_j = f -
\mathop{\sum}_{j=1}^n 0\\
& = f,\end{split}\] and hence the equality follows. This implies
that $\cll$ is a closed subspace and \[P_\cll =
\mathop{\sum}_{i=1}^n O_i = I - \mathop{\Pi}_{i=1}^n (I -
P_i).\]This completes the proof of the lemma. \qed

\newsection{Quotient Modules}

In this section we prove the central result of this paper that a
doubly commuting quotient module of $H^2(\mathbb{D}^n)$ is precisely
the Hilbert tensor product of $n$ quotient modules of the Hardy
module $H^2(\mathbb{D})$.

We begin by generalizing the fact that a closed subspace $\clm$ of
$H^2(\mathbb{D}^n)$ is $M_{z_1}$-reducing if and only if
\[\clm = H^2(\mathbb{D}) \otimes \cle,\]for some closed subspace
$\cle \subseteq H^2(\mathbb{D}^{n-1})$.

\begin{Proposition}\label{1-reducing}
Let $\clq_1$ be a quotient module of $H^2(\mathbb{D})$ and $\clm$ be
a closed subspace of \[\clq = \clq_1 \otimes
\underbrace{H^2(\mathbb{D}) \otimes \cdots \otimes
H^2(\mathbb{D})}_{(n-1) \rm \; times}  \subseteq H^2(\mathbb{D}^n).
\] Then $\clm$ is a $P_{\clq} M_{z_1}|_{\clq}$-reducing subspace of $\clq$ if
and only if \[\clm = \clq_1 \otimes \cle,\]for some closed subspace
$\cle$ of $H^2(\mathbb{D}^{n-1})$.
\end{Proposition}

\NI\textsf{Proof.} Let  $\clm$ be a $P_{\clq}
M_{z_1}|_{\clq}$-reducing subspace of $\clq$. Then
\begin{equation}\label{M-Q}(P_{\clq} M_{z_1}|_{\clq}) P_{\clm} =
P_{\clm} (P_{\clq} M_{z_1}|_{\clq}),\end{equation}or equivalently,
\[(P_{\clq_1} M_{z}|_{\clq_1} \otimes I_{H^2(\mathbb{D})} \otimes
\cdots \otimes I_{H^2(\mathbb{D})}) P_{\clm} = P_{\clm} (P_{\clq_1}
M_{z}|_{\clq_1} \otimes I_{H^2(\mathbb{D})} \otimes \cdots \otimes
I_{H^2(\mathbb{D})}).\] Now \[I_{\clq} - (P_{\clq} M_{z_1}|_{\clq})
(P_{\clq} M_{z_1}|_{\clq})^* = (P_{\clq_1} P_{\mathbb{C}}|_{\clq_1})
\otimes I_{H^2(\mathbb{D})} \otimes \cdots \otimes
I_{H^2(\mathbb{D})}.\]Further (\ref{M-Q}) yields
\[P_{\clm} ((P_{\clq_1} P_{\mathbb{C}}|_{\clq_1}) \otimes I_{H^2(\mathbb{D})}
\otimes \cdots \otimes I_{H^2(\mathbb{D})}) = ((P_{\clq_1}
P_{\mathbb{C}}|_{\clq_1}) \otimes I_{H^2(\mathbb{D})} \otimes \cdots
\otimes I_{H^2(\mathbb{D})}) P_{\clm},\]and therefore
\[P_{\clm} ((P_{\clq_1} P_{\mathbb{C}}|_{\clq_1}) \otimes
I_{H^2(\mathbb{D})} \otimes \cdots \otimes I_{H^2(\mathbb{D})})\] is
the orthogonal projection onto \[\cll : = \clm \cap \mbox{ran~}
((P_{\clq_1} P_{\mathbb{C}}|_{\clq_1}) \otimes I_{H^2(\mathbb{D})}
\otimes \cdots \otimes I_{H^2(\mathbb{D})}) = \clm \cap (\cll_1
\otimes H^2(\mathbb{D}) \otimes \cdots \otimes
H^2(\mathbb{D})),\]where \[\cll_1 = \mbox{ran~} (P_{\clq_1}
P_{\mathbb{C}}|_{\clq_1}) \subseteq \clq_1.\]Since $\cll \subseteq
\cll_1 \otimes H^2(\mathbb{D}) \otimes \cdots \otimes
H^2(\mathbb{D})$ and $\mbox{dim} \cll_1 = 1$ (otherwise, by Lemma
\ref{Q-L} that $\cll_1 = \{0\}$ is equivalent to $\clq_1 = \{0\}$)
we obtain
\[\cll = \cll_1 \otimes \cle,\]for some closed subspace $\cle
\subseteq H^2(\mathbb{D}^{n-1})$. More precisely \[P_{\clm}
((P_{\clq_1} P_{\mathbb{C}}|_{\clq_1}) \otimes I_{H^2(\mathbb{D})}
\otimes \cdots \otimes I_{H^2(\mathbb{D})}) = P_{\cll} = P_{\cll_1
\otimes \cle}.\]We claim that \[\clm = \mathop{\vee}_{l =
0}^{\infty} P_{\clq} M_{z_1}^l \cll.\]Since $\clm$ is $P_{\clq}
M_{z_1}|_{\clq}$-reducing subspace and $\clm \supseteq \cll$, it
follows that
\[\clm \supseteq \mathop{\vee}_{l = 0}^{\infty} P_{\clq} M_{z_1}^l \cll.\] To prove
the reverse inclusion, we let $f \in \clm$ and  $f = \sum_{\bm{k}
\in \mathbb{N}^n} a_{\bm{k}} \bm{z}^{\bm{k}}$, where $a_{\bm{k}} \in
\mathbb{C}$ for all $\bm{k} \in \mathbb{N}^n$. Then
\[ f = P_{\clm} P_{\clq} f = P_{\clm} P_{\clq} \sum_{\bm{k} \in
\mathbb{N}^n} a_{\bm{k}} \bm{z}^{\bm{k}}.\]Observe now that for all
$\bm{k} \in \mathbb{N}^n$,
\[\begin{split} P_{\clm} P_{\clq} \bm{z}^{\bm{k}} & = P_{\clm} ((P_{\clq_1} z_1^{k_1}) (z_2^{k_2}
\cdots z_n^{k_n}))\\ & = P_{\clm} ((P_{\clq_1} M_{z_1}^{k_1}
P_{\clq_1} 1) (z_2^{k_2} \cdots z_n^{k_n}))\\& = (P_{\clm}
P_{\clq} M_{z_1}^{k_1} P_{\clq}) (P_{\clq_1} 1 \otimes z_2^{k_2} \cdots z_n^{k_n}) \\
& = P_{\clq} M_{z_1}^{k_1} (P_{\clm} (P_{\clq_1} 1 \otimes z_2^{k_2}
\cdots z_n^{k_n})),\end{split}\]by (\ref{M-Q}), where the second
equality follows from \[\langle z_1^{k_1}, f \rangle = \langle 1,
(M_{z_1}^{k_1})^* f \rangle = \langle P_{\clq_1} M_{z_1}^{k_1}
P_{\clq_1} 1, f \rangle,\]for all $f \in \clq_1$. Now by applying
Corollary \ref{Q_1}, we obtain that $P_{\clq_1} 1 \in \cll_1$ and
hence \[P_{\clm} (P_{\clq_1} 1 \otimes z_2^{k_2} \cdots z_n^{k_n})
\in \cll.\] Therefore, we infer
\[P_{\clm} P_{\clq} \bm{z}^{\bm{k}} \in \mathop{\vee}_{l = 0}^{\infty} P_{\clq_1}
M_{z_1}^l \cll,\]for all $\bm{k} \in \mathbb{N}^n$ and hence $f \in
\vee_{l = 0}^{\infty} P_{\clq} M_{z_1}^l \cll$. Thus we get $\clm =
\vee_{l=0}^\infty P_{\clq} M_{z_1}^l \cll$. Finally, $\cll = \cll_1
\otimes \cle$ yields \[\clm = \mathop{\vee}_{l = 0}^{\infty}
P_{\clq} M_{z_1}^l \cll = (\mathop{\vee}_{l=0}^{\infty} P_{\clq_1}
M_{z_1}^l \cll_1) \otimes \cle,\]and therefore by Lemma \ref{Q-L},
\[\clm = \clq_1 \otimes \cle.\]The converse part is
trivial. This finishes the proof.  \qed

We are now ready to prove the main result of this section.

\begin{Theorem}\label{dc-Q}
Let $\clq$ be a quotient module of $H^2(\mathbb{D}^n)$. Then $\clq$
is doubly commuting if and only if there exists quotient modules
$\clq_1, \ldots, \clq_n$ of $H^2(\mathbb{D})$ such that \[\clq =
\clq_1 \otimes \cdots \otimes \clq_n.\]
\end{Theorem}

\NI\textsf{Proof.} Let $\clq$ be a doubly commuting quotient module
of $H^2(\mathbb{D}^n)$. Define \[\tilde{\clq}_1 =
\overline{\mbox{span}} \{z_2^{l_2} \cdots z_n^{l_n} \clq : l_2,
\ldots, l_n \in \mathbb{N}\}.\]Then $\tilde{\clq}_1$ is a joint
$(M_{z_2}, \ldots, M_{z_n})$-reducing subspace of
$H^2(\mathbb{D}^n)$. But Proposition \ref{reducing} now allows us to
conclude that
\[\tilde{\clq}_1 = \clq_1 \otimes \underbrace{H^2(\mathbb{D}) \otimes \cdots \otimes
H^2(\mathbb{D})}_{(n-1) \rm\,times},\]for some closed subspace
$\clq_1$ of $H^2(\mathbb{D})$. Since $\tilde{\clq}_1$ is
$M_{z_1}^*$-invariant subspace, that $\clq_1$ is a $M_z^*$-invariant
subspace of $H^2(\mathbb{D})$, that is, $\clq_1$ is a quotient
module of $H^2(\mathbb{D})$.

\NI Note that $\clq \subseteq \tilde{\clq}_1$. We claim that $\clq$
is a $M_{z_1}^*|_{\tilde{\clq}_1}$-reducing subspace of
$\tilde{\clq}_1$, that is, \[P_{\clq} (M_{z_1}^*|_{\tilde{\clq}_1})
= (M_{z_1}^*|_{\tilde{\clq}_1}) P_{\clq}.\]

\NI In order to prove the claim we first observe that for all $l
\geq 0$ and $2 \leq i \leq n$,
\[C_{z_1}^* C_{z_i}^l = C_{z_i}^l C_{z_1}^*,\] and hence
\[C_{z_1}^* C_{z_2}^{l_2} \cdots C_{z_n}^{l_n} = C_{z_2}^{l_2} \cdots C_{z_n}^{l_n}
C_{z_1}^*,\]for all $l_2, \ldots, l_n \geq 0$, that is,
\[M_{z_1}^* P_{\clq} M_{z_2}^{l_2} \cdots M_{z_n}^{l_n} P_{\clq} =
P_{\clq} M_{z_2}^{l_2} \cdots M_{z_n}^{l_n} M_{z_1}^* P_{\clq},\]
or, \[ M_{z_1}^* P_{\clq} M_{z_2}^{l_2} \cdots M_{z_n}^{l_n}
P_{\clq} = P_{\clq} M_{z_1}^* M_{z_2}^{l_2} \cdots M_{z_n}^{l_n}
P_{\clq}.\] From this it follows that for all $f \in \clq$ and $l_2,
\ldots, l_n \geq 0$,\[\begin{split}(P_{\clq}
M_{z_1}^*|_{\tilde{\clq}_1}) (z_2^{l_2} \cdots z_n^{l_n} f) & =
P_{\clq} M_{z_1}^* (z_2^{l_2} \cdots z_n^{l_n} f) = (P_{\clq}
M_{z_1}^* M_{z_2}^{l_2} \cdots M_{z_n}^{l_n}) f \\& = (P_{\clq}
M_{z_1}^* M_{z_2}^{l_2} \cdots M_{z_n}^{l_n} P_{\clq}) f =
M_{z_1}^* P_{\clq} M_{z_2}^{l_2} \cdots M_{z_n}^{l_n} P_{\clq} f \\
& = M_{z_1}^* P_{\clq} M_{z_2}^{l_2} \cdots M_{z_n}^{l_n} f =
(M_{z_1}^* P_{\clq}) (z_2^{l_2} \cdots z_n^{l_n} f).\end{split}\]

\NI Also by $P_{\clq} \tilde{\clq}_1 \subseteq \tilde{\clq}_1$ we
have
\[P_{\clq}
P_{\tilde{\clq}_1} = P_{\tilde{\clq}_1} P_{\clq}
P_{\tilde{\clq}_1}.\] This yields
\[\begin{split}(P_{\clq} M_{z_1}^*|_{\tilde{\clq}_1}) (z_2^{l_2} \cdots
z_n^{l_n} f) & = (M_{z_1}^* P_{\clq}) (z_2^{l_2} \cdots z_n^{l_n} f)
\\ & = M_{z_1}^* P_{\clq} P_{\tilde{\clq}_1} (z_2^{l_2} \cdots
z_n^{l_n} f) \\ & = M_{z_1}^* P_{\tilde{\clq}_1} P_{\clq}
P_{\tilde{\clq}_1} (z_2^{l_2} \cdots z_n^{l_n} f) \\ & =
(M_{z_1}^*|_{\tilde{\clq}_1} P_{\clq}) (z_2^{l_2} \cdots z_n^{l_n}
f),\end{split}\]for all $f \in \clq$ and $l_2, \ldots, l_n \geq 0$,
and therefore
\[P_{\clq} (M_{z_1}^*|_{\tilde{\clq}_1}) = (M_{z_1}^*|_{\tilde{\clq}_1})
P_{\clq}.\]Hence $\clq$ is a $M_{z_1}^*|_{\tilde{\clq}_1}$-reducing
subspace of $\tilde{\clq}_1 = \clq_1 \otimes H^2(\mathbb{D}) \otimes
\cdots \otimes H^2(\mathbb{D})$. Applying Proposition
\ref{1-reducing}, we obtain a closed subspace $\cle_1$ of
$H^2(\mathbb{D}^{n-1})$ such that
\[\clq = \clq_1 \otimes \cle_1.\]
Moreover, since\[\mathop{\vee}_{l = 0}^{\infty} z_1^l \clq =
\mathop{\vee}_{l=0}^{\infty} z_1^l (\clq_1 \otimes \cle_1) =
H^2(\mathbb{D}) \otimes \cle_1,\]and $\mathop{\vee}_{l = 0}^{\infty}
z_1^l \clq $ is a doubly commuting quotient module of
$H^2(\mathbb{D}^n)$, we have that $\cle_1 \subseteq
H^2(\mathbb{D}^{n-1})$ a doubly commuting quotient module of
$H^2(\mathbb{D}^{n-1})$.

\NI By the same argument as above, we conclude that
\[\cle_1 = \clq_2 \otimes \cle_2,\]for some doubly commuting
quotient module of $H^2(\mathbb{D}^{n-2})$. Continuing this process,
we have
\[\clq = \clq_1 \otimes \cdots \otimes \clq_n,\]where $\clq_1,
\ldots, \clq_n$ are quotient modules of $H^2(\mathbb{D})$.

\NI The converse implication follows from the fact that the module
multiplication operators on $\clq = \clq_1 \otimes \cdots \otimes
\clq_n$ are given by \[\{I_{\clq_1} \otimes \cdots \otimes
\underbrace{P_{\clq_i} M_{z}|_{\clq_i}}_{i^{\rm\, th} \rm\,place}
\otimes \cdots \otimes I_{\clq_n}\}_{i=1}^n,\]which is certainly
doubly commuting. This completes the proof. \qed

As a corollary of the above model, we have the following result.

\begin{Corollary}
Let $\clq$ be a closed subspace of $H^2(\mathbb{D}^n)$. Then $\clq$
is doubly commuting quotient module if and only if there exists
$\{\Theta_i\}_{i=1}^n \subseteq H^{\infty}(\mathbb{D})$ such that
each $\Theta_i$ is either inner or the zero function for all $1 \leq
i \leq n$ and \[\clq = \clq_{\Theta_1} \otimes \cdots \otimes
\clq_{\Theta_n}.\]
\end{Corollary}

\NI\textsf{Proof.} Let $\clq$ be a doubly commuting quotient module
of $H^2(\mathbb{D}^n)$. By Theorem \ref{dc-Q}, we know that
\[\clq = \clq_1 \otimes \cdots \otimes \clq_n,\]where $\clq_1,
\ldots, \clq_n$ are quotient modules of $H^2(\mathbb{D})$. For each
$i \in \{1, \ldots, n\}$, if $\clq_i \varsubsetneq H^2(\mathbb{D})$
then
\[\clq_i = \clq_{\Theta_i} =
H^2(\mathbb{D})/\Theta_iH^2(\mathbb{D}),\] for some inner function
$\Theta_i \in H^{\infty}(\mathbb{D})$. Otherwise, $\clq_i =
H^2(\mathbb{D})$ and we define $\Theta_i \equiv 0$ on $\mathbb{D}$
so that \[\clq_i = H^2(\mathbb{D}) =  \clq_{\Theta_i} =
H^2(\mathbb{D})/(0 \cdot H^2(\mathbb{D})).\]The converse part again
follows from Theorem \ref{dc-Q}, and the corollary is proved. \qed

This result was obtained by Izuchi, Nakazi and Seto in \cite{INS}
for the base case $n = 2$ (also see \cite{IN}).

We conclude this section by recording the uniqueness of the tensor
product representations of the doubly commuting quotient modules in
Theorem \ref{dc-Q}.  The same conclusion holds for a more general
framework. Here, we provide a proof using the Hardy space method.

\NI Let $\clq$ be a doubly commuting quotient module of
$H^2(\mathbb{D}^n)$ and \[\clq = \clq_1 \otimes \cdots \otimes
\clq_n = \clr_1 \otimes \ldots \otimes \clr_n,\]for quotient modules
$\{\clq_i\}_{i=1}^n$ and $\{\clr_i\}_{i=1}^n$ of $H^2(\mathbb{D})$.
We claim that $\clq_i = \clr_i$ for all $i$. In fact,
\[\tilde{\clq}_1 := \mathop{\vee}_{l_2, \ldots, l_n \geq 0}^{\infty} z_2^{l_2} \cdots
z_n^{l_n} \clq = \clq_1 \otimes H^2(\mathbb{D}) \otimes \cdots
\otimes H^2(\mathbb{D}) = \clr_1 \otimes H^2(\mathbb{D}) \otimes
\cdots \otimes H^2(\mathbb{D}),\]and \[\mathop{\bigcap}_{i=2}^n
\mbox{ker} M_{z_i}^*|_{\tilde{\clq}_1} = \clq_1 \otimes \mathbb{C}
\otimes \cdots \otimes \mathbb{C} = \clr_1 \otimes \mathbb{C}
\otimes \cdots \otimes \mathbb{C}.\]Consequently, $\clq_1 = \clr_1$
and similarly, for all other $i = 2,\ldots, n$.

\newsection{submodules}

In this section we relate the Hilbert tensor product structure of
the doubly commuting quotient modules to the Beurling like
representations of the corresponding co-doubly commuting submodules.

To proceed further we introduce one more piece of notation.

\NI\textit{\textsf{Let $\Theta_i \in H^{\infty}(\mathbb{D})$ be a
given function indexed by $i \in \{1, \ldots, n\}$. In what follows
by $\tilde{\Theta}_i \in H^{\infty}(\mathbb{D}^n)$ we denote the
extension function defined by
\[\tilde{\Theta}_i(\bm{z}) = \Theta_i(z_i),\]for all $\bm{z} \in
\mathbb{D}^n$.}}

The following provides an explicit correspondence between the doubly
commuting quotient modules and the co-doubly commuting submodules of
$H^2(\mathbb{D}^n)$.

\begin{Theorem}
Let $\clq$ be a quotient module of $H^2(\mathbb{D}^n)$ and $\clq
\neq H^2(\mathbb{D}^n)$. Then $\clq$ is doubly commuting if and only
if there exists inner functions ${\Theta}_{i_j} \in
H^{\infty}(\mathbb{D})$ for $1 \leq i_1 < \ldots < i_m \leq n$ for
some integer $m \in \{1, \ldots, n\}$ such that
\[\clq = H^2(\mathbb{D}^n)/ [\tilde{\Theta}_{i_1} H^2(\mathbb{D}^n) + \cdots +
\tilde{\Theta}_{i_m} H^2(\mathbb{D}^n)],\]where
$\tilde{\Theta}_{i_j}(\bm{z}) = \Theta_{i_j}(z_{i_j})$ for all
$\bm{z} \in \mathbb{D}^n$.
\end{Theorem}

\NI\textsf{Proof.} Let $\clq$ be a doubly commuting quotient module
of $H^2(\mathbb{D}^n)$. Then by Theorem \ref{dc-Q} we have
\[\clq = \clq_{1} \otimes \cdots \otimes \clq_{n},\]where for each
$1 \leq i \leq n$, $\clq_i$ is a submodule of $H^2(\mathbb{D})$.
Choose $1 \leq m \leq n$ such that \[\clq_{i_j} \neq
H^2(\mathbb{D}),\]for $1 \leq i_1 < \ldots < i_m \leq n$. Then
\[\clq = H^2(\mathbb{D}) \otimes \cdots \otimes \clq_{i_1} \otimes
\cdots \otimes \clq_{i_m} \otimes \cdots \otimes
H^2(\mathbb{D}),\]where $\clq_{i_j} \varsubsetneq H^2(\mathbb{D})$
for all $1 \leq i_1 < \ldots < i_m \leq n$. Let
\[\clq_{i_j} = \clq_{{\Theta}_{i_j}} = (\mbox{ran}
M_{\Theta_{i_j}})^{\perp} = \mbox{ran} (I_{H^2(\mathbb{D})} -
M_{{\Theta}_{i_j}} M_{{\Theta}_{i_j}}^*),\]for some inner function
${\Theta}_{i_j} \in H^{\infty}(\mathbb{D})$ for all $j = 1, \ldots,
m$. Let $\tilde{\Theta}_{i_j}$ be the extension of $\Theta_{i_j}$ to
$H^{\infty}(\mathbb{D}^n)$, that is, as a multiplier,
\[M_{\tilde{\Theta}_{i_j}} = I_{H^2(\mathbb{D})} \otimes \cdots
\otimes I_{H^2(\mathbb{D})} \otimes M_{\Theta_{i_j}} \otimes
I_{H^2(\mathbb{D})} \otimes \cdots \otimes
I_{H^2(\mathbb{D})}.\]Hence,
\[I_{H^2(\mathbb{D}^n)} - M_{\tilde{\Theta}_{i_j}}
M_{\tilde{\Theta}_{i_j}}^* = I_{H^2(\mathbb{D})} \otimes \cdots
\otimes I_{H^2(\mathbb{D})} \otimes (I_{H^2(\mathbb{D})} -
M_{{\Theta}_{i_j}} M_{{\Theta}_{i_j}}^*) \otimes I_{H^2(\mathbb{D})}
\otimes \cdots \otimes I_{H^2(\mathbb{D})},\]so that
\[\begin{split}& \mathop{\Pi}_{1 \leq i_1 < \ldots < i_m \leq n}
(I_{H^2(\mathbb{D}^n)} - M_{\tilde{\Theta}_{i_j}} M_{\tilde{\Theta}_{i_j}}^*) \\
&= I_{H^2(\mathbb{D})} \otimes \cdots \otimes (I_{H^2(\mathbb{D})} -
M_{{\Theta}_{i_1}} M_{{\Theta}_{i_1}}^*) \otimes \cdots \otimes
(I_{H^2(\mathbb{D})} - M_{{\Theta}_{i_m}} M_{{\Theta}_{i_m}}^*)
\otimes I_{H^2(\mathbb{D})} \otimes \cdots \otimes
I_{H^2(\mathbb{D})},
\end{split}\]Taking into account that $\clq$ is the range of the
right hand side operator, that is,
\[\begin{split} \clq & =
H^2(\mathbb{D}) \otimes \cdots \otimes \clq_{i_1} \otimes \cdots
\otimes \clq_{i_m} \otimes \cdots \otimes H^2(\mathbb{D})\\& =
\mbox{ran}[ \mathop{\Pi}_{1 \leq i_1 < \ldots < i_m \leq n}
(I_{H^2(\mathbb{D}^n)} - M_{\tilde{\Theta}_{i_j}}
M_{\tilde{\Theta}_{i_j}}^*)],\end{split}\]we deduce readily that
\[P_{\clq^{\perp}} = I_{H^2(\mathbb{D}^n)} - \mathop{\Pi}_{1 \leq
i_1 < \ldots < i_m \leq n} (I_{H^2(\mathbb{D}^n)} -
M_{\tilde{\Theta}_{i_j}} M_{\tilde{\Theta}_{i_j}}^*).\]Consequently,
by Lemma \ref{P-F} we have
\[\clq^{\perp} = \tilde{\Theta}_{i_1} H^2(\mathbb{D}^n) + \cdots +
\tilde{\Theta}_{i_m} H^2(\mathbb{D}^n),\]or \[\clq =
H^2(\mathbb{D}^n)/ [\tilde{\Theta}_{i_1} H^2(\mathbb{D}^n) + \cdots
+ \tilde{\Theta}_{i_m} H^2(\mathbb{D}^n)].\] Conversely, let \[\clq
= H^2(\mathbb{D}^n)/ [\tilde{\Theta}_{i_1} H^2(\mathbb{D}^n) +
\cdots + \tilde{\Theta}_{i_m} H^2(\mathbb{D}^n)].\] Then
\[P_{\clq} = \mathop{\Pi}_{1 \leq i_1 < \ldots < i_m \leq n}
(I_{H^2(\mathbb{D}^n)} - M_{\tilde{\Theta}_{i_j}}
M_{\tilde{\Theta}_{i_j}}^*).\] Also for all $s \neq t$,
\[\begin{split} P_{\clq} M_{z_s} M_{z_t}^* P_{\clq} & = \mathop{\Pi}_{1 \leq
i_1 < \ldots < i_m \leq n} (I_{H^2(\mathbb{D}^n)} -
M_{\tilde{\Theta}_{i_j}} M_{\tilde{\Theta}_{i_j}}^*) M_{z_t}^*
M_{z_s} \mathop{\Pi}_{1 \leq i_1 < \ldots
< i_m \leq n} (I_{H^2(\mathbb{D}^n)} - M_{\tilde{\Theta}_{i_j}} M_{\tilde{\Theta}_{i_j}}^*)\\
& =  P_{\clq} M_{z_t}^* [\mathop{\Pi}_{\substack{ 1 \leq i_1 <
\ldots < i_m \leq n\\i_j \neq t}} (I_{H^2(\mathbb{D}^n)} -
M_{\tilde{\Theta}_{i_j}} M_{\tilde{\Theta}_{i_j}}^*)] \\&
\;\;\;\;\;\quad [\mathop{\Pi}_{\substack{ 1 \leq i_1 < \ldots < i_m
\leq n\\i_j \neq s}} (I_{H^2(\mathbb{D}^n)} -
M_{\tilde{\Theta}_{i_j}} M_{\tilde{\Theta}_{i_j}}^*)] M_{z_s}
P_{\clq}\\& = P_{\clq} M_{z_t}^* \mathop{\Pi}_{1 \leq i_1 < \ldots <
i_m \leq n} (I_{H^2(\mathbb{D}^n)} - M_{\tilde{\Theta}_{i_j}}
M_{\tilde{\Theta}_{i_j}}^*) M_{z_s} P_{\clq} \\& = P_{\clq}
M_{z_t}^* P_{\clq} M_{z_s} P_{\clq}.
\end{split}\]Consequently, for all $s \neq t$ \[C_{z_s} C^*_{z_t} = P_{\clq}
M_{z_s} M_{z_t}^*|_{\clq} = P_{\clq} M_{z_t}^* P_{\clq}
M_{z_s}|_{\clq} = C^*_{z_t} C_{z_s},\]and hence $\clq$ is doubly
commuting. This concludes the proof. \qed

This result is a generalization of Theorem 3.1 of \cite{INS} by
Izuchi, Nakazi and Seto on the base case $n = 2$.

To complete this section, we present the following result concerning
the orthogonal projection formulae of the co-doubly commuting
submodules and the doubly commuting quotient modules of
$H^2(\mathbb{D}^n)$.

\begin{Corollary}\label{S-proj}
Let $\clq$ be a doubly commuting submodule of $H^2(\mathbb{D}^n)$.
Then there exists an integer $m \in \{1, \ldots, n\}$ and inner
functions $\Theta_{i_j} \in H^{\infty}(\mathbb{D})$ such that
\[\clq^{\perp} = \mathop{\sum}_{1 \leq i_1 < \ldots < i_m \leq n}
\tilde{\Theta}_{i_j} H^2(\mathbb{D}^n),\]where
$\tilde{\Theta_i}(\bm{z})= \Theta_{i_j}(z_{i_j})$ for all $\bm{z}
\in \mathbb{D}^n$. Moreover, \[P_{\clq} = I_{H^2(\mathbb{D}^n)} -
\mathop{\Pi}_{j=1}^m (I_{H^2(\mathbb{D}^n)} -
M_{\tilde{\Theta}_{i_j}} M_{\tilde{\Theta}_{i_j}}^*),\]and
\[P_{\clq^\perp} = \mathop{\Pi}_{j=1}^m
(I_{H^2(\mathbb{D}^n)} - M_{\tilde{\Theta}_{i_j}}
M_{\tilde{\Theta}_{i_j}}^*).\]
\end{Corollary}

The above result is the co-doubly commuting submodules analogue of
Beurling's theorem on submodules of $H^2(\mathbb{D})$.

We finally point out that the earlier classifications of the doubly
commuting quotient modules by Izuchi, Nakazi and Seto \cite{INS} has
many deep applications in the study of the submodules and the
quotient modules of the Hardy module over the bidisc (cf.
\cite{III1, III}). Some of these extensions in $n$-variables ($n
\geq 2$) will be addressed in future work. However, the issue of
essential doubly commutativity of the co-doubly commuting submodules
of $H^2(\mathbb{D}^n)$ will be discussed in the forthcoming paper
\cite{JS}.

\vspace{0.2in}

\NI \textsf{Acknowledgement:} We are grateful to the referee for a
very careful reading of the manuscript. In particular, we would like
to thank the referee for improving our proof of Proposition
\ref{1-reducing}.

\end{document}